\documentclass[english,11pt]{article}
\usepackage[latin1]{inputenc}
\usepackage{amsmath,amsthm,amssymb,babel}

\begin{document}
  \newtheorem{teo}{Theorem}
\newtheorem{defin}{Definition}
\newtheorem{prop}{Proposition}
\newtheorem{cor}{Corollary}
\newtheorem{lemma}{Lemma}
\newtheorem{rem}{Remark}
 \def\RR{{\mathbb R} }
\def\CC{{\mathbb C} }
 \newcommand{\huo}{H_0^1(\Omega )}
 \newcommand{\fin}{\nolinebreak\hfill\rule{2mm}{2mm}}
 \newcommand{\mfin}{\nolinebreak\makebox[3cm]{\rule{2mm}{2mm}}}
 \newcommand{\fxr}{$f:X\rightarrow\RR$ }
 \newcommand{\pslc}{({\rm PS})_c }
 \newcommand{\limn}{\displaystyle\lim_{n\rightarrow\infty} }
 \newcommand{\limeps}{\lim_{\varepsilon\searrow 0} }
 \newcommand{\finf}{\displaystyle\underline{f} }
 \newcommand{\fsup}{\displaystyle\overline{f} }
 \newcommand{\intom}{\displaystyle\int_{\Omega} }
 \newcommand{\intm}{\displaystyle\int_{M} }
 \newcommand{\ep}{\varepsilon }
 \newcommand{\intoe}{\displaystyle\int_{\Omega_{\varepsilon} }}
 \newcommand{\nve}{\mid\nabla v_{\varepsilon}\mid^2 }
 \newcommand{\sd}{\displaystyle\sum_{j=1}^d }
 \newcommand{\lue}{\displaystyle\log\, \frac{1}{\ep} }
 \newcommand{\dn}{\displaystyle\frac{\partial}{\partial\nu} }
 \newcommand{\sdo}{\displaystyle\sup_{\partial\omega_j} }
 \newcommand{\ido}{\displaystyle\inf_{\partial\omega_j} }
 \def\proof{{\it Proof.}\ }
\def\pt#1{{\it Proof of Theorem \ref{#1}.}}
\def\eq#1{(\ref{#1})}
\def\th#1{Theorem \ref{#1}}
\def\neweq#1{\begin{equation}\label{#1}}
\def\endeq{\end{equation}}
\def\weak{\rightharpoonup}
\def\ep{\varepsilon}
\def\half{{1\over2}}
\def\phi{\varphi}
\def\intom{\int_\Omega}

\def\proof{{\it Proof.}\ }
\def\pt#1{{\it Proof of Theorem \ref{#1}.}}
\def\eq#1{(\ref{#1})}
\def\th#1{Theorem \ref{#1}}
\def\neweq#1{\begin{equation}\label{#1}}
\def\endeq{\end{equation}}
\def\weak{\rightharpoonup}
\def\ep{\varepsilon}
\def\half{{1\over2}}
\def\phi{\varphi}
\def\intom{\int_\Omega}
 
 \newcommand{\intoj}{\displaystyle\int_{\partial\omega_j} }
 \newcommand{\di}{\displaystyle} \newcommand{\ri}{\rightarrow}
 \newcommand{\sij}{\sum_{i,j}} \font\rmg=cmr8 \font\bfg=cmbx12
 \newcommand{\incep}{\left\{\begin{array}{cl} }
 \newcommand{\termin}{\end{array}\right. }
 \newcommand{\sfp}{\displaystyle (f(a+)-f(a)) }
 \newcommand{\sfm}{\di (f(a)-f(a-)) }
 \newcommand{\intomp}{\di \int_{\Omega_a^+} }
 \newcommand{\intomm}{\di \int_{\Omega_a^-} }
 \newcommand{\alin}{\begin{array}{cl} }
 \def\dpx{{\rm div}(|\nabla u|^{p(x)-2}\nabla u)}
\font\rmg=cmr10
\font\bfg=cmbx10

\setcounter{page}{1}

\title{\sc On a nonlinear eigenvalue problem in Sobolev spaces with variable exponent}
\author{Teodora-Liliana Dinu \\ \small Department of Mathematics, ``Fra\c tii Buze\c sti" College, 
Bd. \c Stirbei--Vod\u a No. 5, 200352 Craiova, Romania\\ \small 
E-mail: {\tt tldinu@gmail.com}}
  \date{}
      \maketitle

\bigskip
 {\bf Abstract.} We consider a class of nonlinear Dirichlet problems involving
 the $p(x)$--Laplace operator. Our framework is based on the theory of Sobolev spaces with 
variable exponent and we establish the existence of a weak solution in such a space. The proof
relies on the Mountain Pass Theorem.



\normalsize

\section{Introduction}
The Mountain Pass Theorem is due to Ambrosetti and Rabinowitz \cite{1} and is one of the most 
powerful tools in Nonlinear Analysis for proving the existence of critical points of energy 
functionals. One of the simplest versions of the Mountain Pass Theorem asserts that if a 
continuously differential functional has two local minima, then (under some natural assumptions) 
such a function has a third critical point. This fact is elementary for functions of one real 
variable. However, even for functions on the plane the proof of such a theorem requires deep 
topological ideas. The Mountain Pass Theorem has numerous generalizations and has been applied in 
the treatment of various classes of boundary value problems. We refer to the recent monograph by 
Jabri \cite{jabri} for an excellent survey of some of the most interesting applications of this 
abstract result. We do not intend to insist on the wide spectrum of applications of the Mountain 
Pass Theorem. We remark only that this theorem has been applied in the last few years 
in very concrete situations. For instance, in Lewin \cite{lew} it is considered 
 a neutral molecule that possesses two distinct stable positions for its nuclei, and it is looked 
for a mountain pass point between the two minima in the non-relativistic Schr\"odinger framework.

As showed in \cite{1}, one of the simplest applications of the Mountain Pass Theorem implies the 
existence of solutions for the Dirichlet problem
$$\left\{\begin{array}{lll}
&\di -\Delta u=u^{p-1}&\di\qquad\mbox{in $\Omega$}\\
&\di u>0&\di\qquad\mbox{in $\Omega$}\\
&\di u=0&\di\qquad\mbox{on $\partial\Omega$,}
\end{array}\right.
$$
where $\Omega\subset\RR^N$ is a smooth bounded domain, $2<p<2N/(N-2)$ if $N\geq 3$ and $p\in 
(2,\infty)$ if $N=1$ or $N=2$. 
Under the same assumptions on $p$, similar arguments show that the boundary value problem
$$\left\{\begin{array}{lll}
&\di -\Delta u-\lambda u=u^{p-1}&\di\qquad\mbox{in $\Omega$}\\
&\di u>0&\di\qquad\mbox{in $\Omega$}\\
&\di u=0&\di\qquad\mbox{on $\partial\Omega$,}
\end{array}\right.$$
has a solution  for any $\lambda<\lambda_1$, where $\lambda_1$ denotes the first eigenvalue of 
$(-\Delta)$ in $H^1_0(\Omega)$. 
The proof of this result relies on the fact that the operator $(-\Delta -\lambda I)$ is coercive 
if $\lambda<\lambda_1$.
Moreover, a multiplication by $\varphi_1$ and integration on $\Omega$ implies that there is no 
solution if $\lambda\geq\lambda_1$,  where $\varphi_1$ stands for the first eigenfunction of the 
Laplace operator. We refer to \cite{pre} for interesting localization results of solutions to 
problems of the above type, as well as for a lower bound of all nontrivial solutions.

The main purpose of this paper is to study a related problem, but for a more general differential 
operator, the so-called $p(x)$--Laplace operator. This degenerate differential operator is defined 
by $\Delta_{p(x)}u:={\rm div}(|\nabla u|^{p(x)-2}\nabla u)$ (where $p(x)$ is a certain function 
whose properties will be stated in what follows) and that generalizes the celebrated $p$--Laplace 
operator
$\Delta_p u:=\mbox{div} (|\nabla u|^{p-2}\nabla u)$, where $p>1$ is a constant.
The $p(x)$--Laplace operator possesses more complicated nonlinearity than the $p$--Laplacian, for 
example, it is inhomogeneous.
We only recall that $\Delta_p$ describes a  variety of phenomena in the nature. For instance, the 
equation
governing the motion of a fluid involves the $p$--Laplace operator. More
exactly, the shear
stress $\vec\tau$ and the velocity gradient $\nabla u$ of the fluid
are related in the manner that
$\vec\tau (x)=r(x)|\nabla u|^{p-2}\nabla u$,
where $p=2$ (resp., $p<2$ or $p>2$) if the fluid is Newtonian (resp.,
pseudoplastic or dilatant). Other applications of the $p$--Laplacian also appear
 in the study of flow through porous media ($p=3/2$),
Nonlinear Elasticity ($p\geq 2$), or Glaciology
 ($1<p\leq 4/3$). 

\section{Auxiliary results}
In this section we recall the main properties of Lebesgue and Sobolev spaces with variable 
exponent.
  We point out that these functional spaces appeared in the literature for the first time already 
in a 1931 article by W.~Orlicz \cite{orl}, who proved various results (including H\"older's 
inequality) in a discrete framework.
Orlicz also considered the variable exponent function space $L^{p(x)}$ on the real line, and 
proved the H\"older inequality in this setting, too. Next, Orlicz abandoned the study of variable 
exponent spaces, to concentrate on the theory of the function spaces that now bear his name.
The first systematic study of spaces with variable exponent (called {\it modular spaces}) is due 
to Nakano \cite{nak}. In the appendix of this book, Nakano mentions explicitly variable exponent 
Lebesgue spaces as an example of the more general spaces he considers \cite[p.~284]{nak}.
Despite their broad interest, these spaces have not reached the same main-stream position as 
Orlicz spaces. Somewhat later, a more explicit version of such spaces, namely modular function 
spaces, were investigated by Polish mathematicians. We refer to the book by Musielak \cite{M} for 
a nice presentation of modular function spaces. This book, although not dealing specifically with 
the spaces that interest us, is still specific enough to contain several interesting results 
regarding variable exponent spaces. 
Variable exponent Lebesgue spaces on the real line have been independently developed by Russian 
researchers, notably Sharapudinov. These investigations originated in a paper by Tsenov   
\cite{tse}. The question raised by Tsenov and solved by Sharapudinov \cite{sha} is the 
minimization of
$\int_a^b \vert u(x)-v(x)\vert^{p(x)}dx, $
where $u$ is a fixed function and $v$ varies over a finite dimensional subspace of 
$L^{p(x)}([a,b])$.  Sharapudinov also introduces the Luxemburg norm for the Lebesgue space and 
shows that this space is reflexive if the exponent satisfies $1 < p^- \leq p^+ < \infty$. In the 
80's Zhikov started a new line of investigation, that was to become intimately related to the 
study of variable exponent spaces, namely he considered variational integrals with non-standard 
growth conditions. 

Let $\Omega$ be
a bounded open set in $\RR^N$. 

Set
$$C_+(\overline\Omega)=\{h;\;h\in C(\overline\Omega),\;h(x)>1\;{\rm for}\;
{\rm all}\;x\in\overline\Omega\}.$$
For any $h\in C_+(\overline\Omega)$ we define
$$h^+=\sup_{x\in\Omega}h(x)\qquad\mbox{and}\qquad h^-=
\inf_{x\in\Omega}h(x).$$
For any $p(x)\in C_+(\overline\Omega)$, we define the variable exponent Lebesgue and Sobolev 
spaces
$$L^{p(x)}(\Omega)=\{u;\ u\ \mbox{is a
 measurable real-valued function such that }
\int_\Omega|u(x)|^{p(x)}\;dx<\infty\}$$
and
$$W^{1,p(x)}(\Omega)=\{u\in L^{p(x)}(\Omega);\ |\nabla u|\in L^{p(x)}(\Omega)\}\,.$$

On these spaces we define, respectively, the following norms 
$$|u|_{p(x)}=\inf\left\{\mu>0;\;\int_\Omega\left|
\frac{u(x)}{\mu}\right|^{p(x)}\;dx\leq 1\right\}\qquad\mbox{(called Luxemburg norm)}$$
and
$$\|u\|=|u|_{p(x)}+|\nabla u|_{p(x)}\,.$$    
Variable exponent Lebesgue and Sobolev spaces resemble classical Lebesgue and Sobolev spaces
in many respects: they are Banach spaces \cite[Theorem 2.5]{KR}, the H\"older inequality holds
\cite[Theorem 2.1]{KR}, they are reflexive if and only if $1 < p^-\leq p^+<\infty$  
\cite[Corollary 2.7]{KR} and continuous functions are dense if $p^+ <\infty$ \cite[Theorem 
2.11]{KR}. The inclusion between
Lebesgue spaces also generalizes naturally \cite[Theorem 2.8]{KR}: if $0 < |\Omega|<\infty$
 and $p_1,\, p_2\in C_+(\overline\Omega)$
are variable exponent so that $p_1(x) \leq p_2(x)$ in $\Omega$ then there exists the continuous 
embedding
$L^{p_2(x)}(\Omega)\hookrightarrow L^{p_1(x)}(\Omega)$, whose norm does not exceed $|\Omega|+1$.

We denote by $L^{p'(x)}(\Omega)$ the conjugate space 
of $L^{p(x)}(\Omega)$, where $1/p(x)+1/p'(x)=1$. For any 
$u\in L^{p(x)}(\Omega)$ and $v\in L^{p'(x)}(\Omega)$ the H\"older 
type inequality
\begin{equation}\label{Hol}
\left|\int_\Omega uv\;dx\right|\leq\left(\frac{1}{p^-}+
\frac{1}{(p')^-}\right)|u|_{p(x)}|v|_{p'(x)}
\end{equation}
holds true.  

An important role in manipulating the generalized Lebesgue--Sobolev spaces is played by the {\it 
modular} of the $L^{p(x)}(\Omega)$ space, which is the mapping 
 $\rho_{p(x)}:L^{p(x)}(\Omega)\rightarrow\RR$ defined by
$$\rho_{p(x)}(u)=\int_\Omega|u|^{p(x)}\;dx.$$
If $(u_n)$, $u\in L^{p(x)}(\Omega)$ and $p^+<\infty$ then the following relations 
holds true
$$
|u|_{p(x)}>1\;\;\;\Rightarrow\;\;\;|u|_{p(x)}^{p^-}\leq\rho_{p(x)}(u)
\leq|u|_{p(x)}^{p^+}
$$  
\begin{equation}\label{L5}
|u|_{p(x)}<1\;\;\;\Rightarrow\;\;\;|u|_{p(x)}^{p^+}\leq
\rho_{p(x)}(u)\leq|u|_{p(x)}^{p^-}
\end{equation}
$$
|u_n-u|_{p(x)}\rightarrow 0\;\;\;\Leftrightarrow\;\;\;\rho_{p(x)}
(u_n-u)\rightarrow 0.
$$
Spaces with $p^+ =\infty$ have been studied by Edmunds, Lang and Nekvinda \cite{edm}.

Denote by $W_0^{1,p(x)}(\Omega)$ the closure of 
$C_0^\infty(\Omega)$ in $W^{1,p(x)}(\Omega)$. On this space we can use the equivalent norm 
$\|u\|=|\nabla u|_{p(x)}.$
The space $(W_0^{1,p(x)}(\Omega),\|\cdot\|)$ is a separable and 
reflexive Banach space. The dual of this space is denoted by $W_0^{-1,p'(x)}(\Omega)$.  We note 
that if $q\in C_+(\overline\Omega)$ 
and $q(x)<p^\star(x)$ for all $x\in\overline\Omega$ then the 
embedding
$W_0^{1,p(x)}(\Omega)\hookrightarrow L^{q(x)}(\Omega)$ 
is compact, while $W_0^{1,p(x)}(\Omega)$ is continuously embedded into $L^{p^*(x)}(\Omega)$, 
where $p^\star(x)$ denotes the critical Sobolev exponent, that is,
$p^\star(x)=Np(x)/(N-p(x))$, provided that $p(x)<N$ for all $x\in\overline\Omega$.
We refer to \cite{dien,edm2,edm3,FH,FZZ,FZ1,KR,ruz} for further properties and applications of 
variable exponent Lebesgue--Sobolev spaces.

\section{The main result}
Assume throughout this paper that $\Omega$ is a smooth bounded open set in $\RR^N$ ($N\geq 2$), 
$\lambda$ is a real parameter and $p\in C_+(\overline\Omega)$.

Consider the boundary value problem
\begin{equation}\label{8}
\left\{\begin{array}{lll}
\di & -\dpx  = \lambda u^{p(x)-1}+u^{q-1}&\qquad\mbox{in}\ 
\Omega\\
& u=0&\qquad\mbox{on}\ \partial\Omega\\
& u\geq 0,\ u\not\equiv 0&\qquad\mbox{in}\ \Omega\, ,
\end{array}\right.
\end{equation}
where $p\in C_+(\overline\Omega)$ such that $p^+<N$, and $q$ is a real number.

\begin{defin}\label{def1}
Let $\lambda$ be a real number. We say that $u\in W_0^{1,p(x)}(\Omega)$ is a solution of Problem 
\eq{8} if $u\geq 0$, $u\not\equiv 0$ in $\Omega$ and
$$\intom |\nabla u|^{p(x)-2}\nabla u\nabla vdx=\lambda\intom u^{p(x)-1}vdx+\intom u^{q-1}vdx,
\qquad\forall v\in W_0^{1,p(x)}(\Omega)\,.$$
\end{defin}

A crucial role in the statement of our result will be played by the nonlinear eigenvalue problem
\begin{equation}\label{8pr}
\left\{\begin{array}{lll}
\di & -\dpx  = \lambda |u|^{p(x)-2}u&\qquad\mbox{in}\ 
\Omega\\
& u=0&\qquad\mbox{on}\ \partial\Omega\, .
\end{array}\right.
\end{equation}
It follows easily that if $(u,\lambda)$ is a solution of \eq{8pr} and $u\not\equiv 0$ then
$$\lambda =\lambda (u)=\frac{\intom |\nabla u|^{p(x)}dx}{\intom |u|^{p(x)}dx}$$
and hence $\lambda >0$. Let $\Lambda$ denote the set of eigenvalues of \eq{8pr}, that is,
$$\Lambda =\Lambda_{p(x)}=\{\lambda\in\RR;\ \lambda\ \mbox{is an eigenvalue of Problem 
\eq{8pr}}\}\,.$$
In \cite{garper} it is showed that if the function $p(x)$ is a constant $p>1$ (we refer to 
\cite{hbrezis} for the linear case $p(x)\equiv 2$), then Problem \eq{8pr}
has a sequence of eigenvalues, $\sup\Lambda =+\infty$ and $\inf\Lambda=\lambda_1=\lambda_{1,p}>0$,
where $\lambda_{1,p}$ is the first eigenvalue of $(-\Delta_p)$ in $W_0^{1,p}(\Omega)$ and
$$\lambda_1=\lambda_{1,p}=\inf_{u\in W_0^{1,p}(\Omega\setminus\{0\})}\frac{\intom |\nabla 
u|^{p(x)}dx}{\intom |u|^{p(x)}dx}\,.$$
In \cite{FZZ} it is showed that for general functions $p(x)$ the set $\Lambda$ is infinite 
and $\sup\Lambda=+\infty$. Moreover, it may arise that $\inf\Lambda=0$. Set
$$\lambda^*=\lambda^*_{p(x)}=\inf\Lambda\,.$$
In \cite{FZZ} it is argued that if $N=1$ then $\lambda^*>0$ if and only if the function $p(x)$
is monotone. In arbitrary dimension, $\lambda^*=0$ provided that there exist an open set 
$U\subset\Omega$ and a point $x_0\in U$ such that $p(x_0)<$ (or $>$) $p(x)$ for all $x\in\partial 
U$. 

\begin{teo}\label{t1}  Assume that $\lambda <\lambda^*$ and $p^+<q<Np^-/(N-p^-)$.
Then Problem \eq{8} has at least a 
solution.\end{teo}

We cannot expect that Problem \eq{8} has a solution for any $\lambda\geq\lambda^*$. Indeed, 
consider the simplest case $p(x)\equiv 2$, take $\lambda\geq\lambda_1$ and multiply the equation
in \eq{8} by $\varphi_1>0$. Integrating on $\Omega$ we find
$$(\lambda-\lambda_1)\intom u\varphi_1dx+\intom u^{q-1}\varphi_1dx=0$$
which yields a contradiction.

The proof of the above result relies on the celebrated Mountain Pass Theorem of Ambrosetti and 
Rabinowitz \cite{1} in the following variant.

\begin{teo}\label{ambr}
 Let $X$ be a real Banach space and let
$F:X\ri\RR$ be a $C^1$--functional. Suppose that $F$ satisfies the
Palais-Smale condition and the following geometric assumptions:

\begin{equation}\label{9}\left\{\begin{array}{cl}
& \mbox{there exist positive constants $R$ and $c_0$ such that}\\
&F(u)\geq c_0 ,\ \mbox{for all $u\in X$ with $\| u\| =R$}\, ;
\end{array}\right.
\end{equation}

\begin{equation}\label{10}F(0)<c_0 \ \mbox{and there exists }
 v\in X\ \mbox{such that $\| v\|
>R$ and $F(v)<c_0$}\, .
\end{equation}

Then the functional $F$ possesses at least a critical point.\end{teo}

We recall the celebrated ``compactness condition" introduced by Palais and Smale
\cite{palsma}: the functional $F\in C^1(X,\RR)$ satisfies the Palais-Smale condition provided that 
any sequence $(u_n)$ in $X$
such that $$\sup_n|F(u_n)|<\infty\qquad\mbox{ and }\qquad\| F'(u_n)\|\ri 0$$ has a convergent 
subsequence.

The name of the above result is a consequence of a simplified visualization for the objects from 
theorem. Indeed,  consider the set $ \{ 0, \, v\}$, where $0$ and $v$ are two villages, and the 
set of all paths joining $ 0$ and $ v$. Then, assuming that $ F(u)$ represents the altitude of 
point $u$, assumptions \eq{9} and \eq{10} are equivalent to say that the villages $0$ and $v$ are 
separated by a mountains chain. So, the conclusion of the theorem tells us that there exists a 
path between the villages with a minimal altitude. With other words, there exists a ``mountain 
pass''.

\section{Proof of Theorem \ref{t1}}

Our hypothesis  $\lambda <\lambda^*$
implies that there exists $C_0>0$ such that
\begin{equation}\label{35}
\intom (|\nabla v|^{p(x)}-\lambda| v|^{p(x)})dx\geq C_0\intom|\nabla
v|^{p(x)}dx\qquad\mbox{for all $v\in
W_0^{1,p(x)}(\Omega )$}\, .\end{equation}

Set
$$g(u)=
\incep
& u^{q-1},\quad\mbox{if $u\geq 0$}\, ,\\
&0,\quad\mbox{if $u<0$}
\termin
$$
and $G(u)=\int_0^ug(t)dt$. Define the energy functional associated to Problem \eq{8} by
$$J(u)=\intom \frac{1}{p(x)}\,\left(|\nabla u|^{p(x)}-\lambda| u|^{p(x)}\right)dx-\intom
G(u)dx\qquad\mbox{for all $u\in W_0^{1,p(x)}$}\, .$$
Observe that
$$| G(u)|\leq C\,| u|^q$$
and, by our hypotheses on $p(x)$ and $q$, we have $W_0^{1,p(x)}(\Omega )\hookrightarrow L^q(\Omega
)$, which implies that $J$ is well defined on $W_0^{1,p(x)}(\Omega )$.

A straightforward computation shows that $J$ is of class $C^1$ and, for
every $v\in W_0^{1,p(x)}(\Omega )$,
$$J'(u)(v)=\intom (|\nabla u|^{p(x)-2}\nabla u\cdot\nabla
v-\lambda\,| u|^{p(x)-2}uv)dx-\intom g(u)vdx\, .$$

We prove in what follows that $J$ satisfies the hypotheses of the
 Mountain Pass Theorem.

\medskip
{\sc Verification of \eq{9}}.  
We may write, for every $u\in\RR$,
$$| g(u)|\leq | u|^{q-1}\, .$$
Thus, for every $u\in\RR$,
\begin{equation}\label{36}| G(u)|\leq\frac{1}{q}\, | u|^q\,
.\end{equation}
Next, by \eq{35} and \eq{36},
\begin{equation}\label{37}
 J(u)\di\geq \frac{C_0}{p^+}\,\intom |\nabla u|^{p(x)}dx-C\,\intom |u|^{q}dx
= C_1\,\intom |\nabla u|^{p(x)}dx-C_2\, \| u\|_{L^q}^{q}dx
\,,\end{equation}
for every $u\in W_0^{1,p(x)}(\Omega )$, where $C_1$ and $C_2$ are positive constants.
So, by relation \eq{L5} and 
using the compact embedding $W_0^{1,p(x)}(\Omega)\hookrightarrow L^q(\Omega)$
combined with the assumption $p^+<q$ we find,
for all $u\in W_0^{1,p(x)}(\Omega )$ with $\| u\|=|\nabla u|_{p(x)} =R$ sufficiently small,
$$J(u)\geq C_1\, |\nabla u|_{p(x)}^{p^+}-C_3\,\, |\nabla u|_{p(x)}^{q}\geq c_0>0\,.$$

\medskip
{\sc Verification of \eq{10}}.  Choose $u_0\in W_0^{1,p(x)}(\Omega )$, $u_0>0$ in
$\Omega$. Since $p^+<q$, it follows that if $t>0$ is large enough then
$$J(tu_0)=\intom \frac{t^{p(x)}}{p(x)}\,\left(|\nabla u_0|^{p(x)}-\lambda|
u_0|^{p(x)}\right)dx-\frac{t^q}{q}\,\intom u_0^qdx<0\, .$$

\medskip
{\sc Verification of the Palais-Smale condition}. Let $(u_n)$ be a sequence
in $W_0^{1,p(x)}(\Omega )$ such that
\begin{equation}\label{38}\sup_n| J(u_n)| <+\infty\end{equation}
\begin{equation}\label{39}\| J'(u_n)\|_{W^{-1,p'(x)}}\ri 0\qquad\mbox{as $n\ri\infty$}\, 
.\end{equation}

We first prove that $(u_n)$ is bounded in $W_0^{1,p(x)}(\Omega )$. Remark that
\eq{39} implies that, for every $v\in W_0^{1,p(x)}(\Omega )$,
\begin{equation}\label{40}
\intom (|\nabla u_n|^{p(x)-2}\nabla u_n\cdot\nabla v-\lambda \,|
u_n|^{p(x)-2}\, u_nv)dx
=\intom g(u_n)vdx+o(1)\,\| v\|\, \qquad \mbox{ as
$n\ri\infty$}\, .
\end{equation}
Choosing $v=u_n$ in \eq{40} we find
\begin{equation}\label{41}\intom \left(|\nabla u_n|^{p(x)}-\lambda\, | u_n|^{p(x)}\right)dx=\intom
g(u_n)u_ndx+o(1)\,\| u_n\|\, .\end{equation}
Relation \eq{38} implies that there exists $M>0$ such that, for any $n\geq 1$,
\begin{equation}\label{42}\left|\intom \frac{1}{p(x)}\,\left(|\nabla u_n|^{p(x)}-\lambda\,| 
u_n|^{p(x)}\right)dx-\intom
G(u_n)dx\right|\leq M\, .\end{equation}
But a simple computation yields
\begin{equation}\label{43}\intom g(u_n)u_ndx=q\,\intom G(u_n)dx\, .\end{equation}
Combining \eq{41}, \eq{42} and \eq{43} and using our assumption $p^+<q$ we find
\begin{equation}\label{44}\intom G(u_n)dx=O(1)+o(1)\,\| u_n\|\, .\end{equation}
Thus, by \eq{41} and \eq{44},
$$\intom |\nabla u_n|^{p(x)}dx= O(1)+o(1)\,\| u_n\|\, ,$$
which means that $(u_n)$ is bounded in $W_0^{1,p(x)}(\Omega)$.

It remains to prove that $(u_n)$ is relatively compact. 
We first  remark
that \eq{40} may be rewritten as
\begin{equation}\label{45}\intom |\nabla u_n|^{p(x)-2}\,\nabla u_n\cdot\nabla vdx
=\intom
h(x,u_n)vdx+o(1)\,\| v\|\, ,\end{equation}
for every $v\in W_0^{1,p(x)}(\Omega )$, where
$$h(x,u)=g(u)+\lambda\,| u|^{p(x)-2}\, u\, ,$$
where $\lambda <\lambda^*$ is fixed.
Obviously, $h$ is continuous and, since $q<Np(x)/(N-p(x))$ for all $x\in\overline\Omega$, there 
exists $C>0$ such that
\begin{equation}\label{46}| h(x,u)|\leq C\, \left(1+| u|^{(Np(x)-N+p(x))/(N-p(x))}\, \right)\qquad 
\mbox{for all $x\in\overline\Omega$ and $u\in\RR$}\, .\end{equation}
Moreover
\begin{equation}\label{47}h(x,u)=o\left(| u|^{Np(x)/(N-p(x))}\right)\qquad\mbox{as $| 
u|\ri\infty$,
uniformly for $x\in\overline\Omega$}\,
.\end{equation}

Define $A:W_0^{1,p(x)}(\Omega)\ri W^{-1,p'(x)}(\Omega)$
by $Au=-\dpx$. Then $A$ is invertible and
$A^{-1}:W^{-1,p'(x)}(\Omega )\ri W_0^{1,p(x)}(\Omega )$
is a continuous operator. Thus, by \eq{45}, it suffices to show that
$h(x,u_n)$ is relatively compact in $ W^{-1,p'(x)}(\Omega )$. By continuous embeddings for
Sobolev spaces with variable exponent,
this will be achieved by proving that a subsequence of $h(x,u_n)$ is
convergent in $$(L^{Np(x)/(N-p(x))}(\Omega )\, )^\star =L^{Np(x)/(Np(x)-N+p(x))}(\Omega
)\,.$$

Since $(u_n)$ is bounded in $ W_0^{1,p(x)}(\Omega )\subset
L^{Np(x)/(N-p(x))}(\Omega )$ we can suppose that, up to a subsequence,
$$u_n\ri u\in L^{Np(x)/(N-p(x))}(\Omega ) \qquad\mbox{a.e. in $\Omega$}\, .$$
Moreover, by Egorov's Theorem, for each $\delta >0$, there exists a subset $A$
of $\Omega$ with $| A| <\delta$ and such that
$$u_n\ri u\, \qquad\mbox{uniformly in $\Omega\setminus A$}\, .$$
So, it is sufficient to show that 
$$\int_A| h(u_n)-h(u)|^{Np(x)/(Np(x)-N+p(x))}\, dx\leq\eta\, ,$$
for any fixed $\eta >0$. But, by \eq{46},
$$\int_A| h(u)|^{Np(x)/(Np(x)-N+p(x))}dx\leq C\,\int_A(1+|
u|^{Np(x)/(N-p(x))})dx\, ,$$
which can be made arbitrarily small if we choose a sufficiently small $\delta
>0$.

We have, by \eq{47},
$$\int_A| h(u_n)-h(u)|^{Np(x)/(Np(x)-N+p(x))}dx\leq\ep\,\int_A|
u_n-u|^{Np(x)/(N-p(x))}dx+C_\ep\, | A|\, ,$$
which can be also made arbitrarily small, by continuous embeddings for Sobolev spaces with 
variable exponent combined with the
boundedness of $(u_n)$ in $W_0^{1,p(x)}(\Omega )$.
Hence, $J$ satisfies the Palais-Smale condition. Thus, by Theorem \ref{ambr},
the boundary value problem 
$$
\left\{\begin{array}{lll}
\di & -\dpx  = \lambda |u|^{p(x)-2}u+g(u)&\qquad\mbox{in}\ 
\Omega\\
& u=0&\qquad\mbox{on}\ \partial\Omega
\end{array}\right.$$
 has a weak solution $u\in W_0^{1,p(x)}(\Omega)\setminus\{0\}$. It remains to show that $u\geq 0$. 
Indeed, multiplying the equation by $u^-$ and integrating we find
$$\intom |\nabla u^-|^{p(x)}dx-\lambda\intom (u^-)^{p(x)}dx=0\,.$$
Thus, since $\lambda <\lambda^*$, we deduce that $u^-=0$ in $\Omega$ or, equivalently, $u\geq 0$
in $\Omega$.
\qed

\medskip
A careful analysis of the above proof shows that the existence result stated in Theorem 
\ref{t1} remains valid if $u^{q-1}$ is replaced by the more general nonlinearity $f(x,u)$,
where $f(x,u):\overline\Omega\ri\RR$ is a continuous functions satisfying
$$|f(x,u)|\leq C(|u|+|u|^{q-1}),\qquad\forall x\in\Omega,\ \forall u\in\RR$$
with $p^-<q<Np^+/(N-p^+)$ if $N\geq 3$ and $q\in (p^-,\infty)$ if $N=1$ or $N=2$, 
$$\lim_{\ep\searrow 0}\sup\left\{\left|\frac{f(x,t)}{t}\right|;\ (x,t)\in\overline\Omega\times
(-\ep,\ep)\right\}=0\qquad\mbox{uniformly for $x\in\overline\Omega$}$$
and
$$0\leq\mu F(x,u)\leq uf(x,u)\qquad\mbox{for $0<u$ large and some $\mu>2$},$$
where $F(x,u)=\int_0^uf(x,t)dt$. 

The following result shows that Theorem \ref{t1} still remains valid if the right hand-side is 
affected by a small perturbation. Consider the boundary value problem
\begin{equation}\label{8cor}
\left\{\begin{array}{lll}
\di & -\dpx  = \lambda |u|^{p(x)-2}u+|u|^{q-2}u+a(x)&\qquad\mbox{in}\ 
\Omega\\
& u=0&\qquad\mbox{on}\ \partial\Omega\, ,
\end{array}\right.
\end{equation}
where $a\in L^\infty (\Omega)$, $p\in C_+(\overline\Omega)$ such that $p^+<N$, and $q$ is a real 
number.

\begin{cor}\label{cor111}  Assume that $\lambda <\lambda^*$ and $p^+<q<Np^-/(N-p^-)$.
There exists $\delta>0$ such that if $\|a\|_{L^\infty}<\delta$ then Problem \eq{8cor} has at least 
a solution.\end{cor}

\proof For any $u\in W_0^{1,p(x)}(\Omega)$ define the energy functional
$$E(u)=\intom \frac{1}{p(x)}\,\left(|\nabla u|^{p(x)}-\lambda| u|^{p(x)}\right)dx-\frac 1q\intom
|u|^qdx-\intom a(x)udx\,.$$
We have already seen that if $a=0$ then Problem \eq{8cor} has a nontrivial and nonnegative 
solution. If $\|a\|_{L^\infty}$ is sufficiently small then the verification of the Palais-Smale 
condition, as well as of the geometric assumptions \eq{9} and \eq{10} can be made following the 
same ideas as in the proof of Theorem \ref{t1}. Thus, by Theorem \ref{ambr}, the functional $E$ 
has a nontrivial critical point $u\in W_0^{1,p(x)}(\Omega)$, which is a solution of Problem 
\eq{8cor}. However, we are not able to decide if this solution is nonnegative. This result remains 
true if $a\geq 0$, as we can see easily after multiplication with $u^-$ and integration. \qed


\begin{thebibliography}{99}  {\footnotesize

\bibitem{1} A. Ambrosetti and P. Rabinowitz, Dual variational methods in
critical point theory and applications, {\it J.~Funct. Anal.} {\bf 14} (1973),
349-381.

\bibitem{hbrezis} H. Brezis, {\it Analyse fonctionnelle: th\'eorie et
applications}, Masson, Paris, 1992.

\bibitem{dien} L. Diening, Theorical and numerical results for electrorheological fluids, Ph.D. 
thesis, University of Freiburg, Germany, 2002.

\bibitem{edm} D. E. Edmunds, J. Lang, and A. Nekvinda, On $L^{p(x)}$ norms, {\it Proc. Roy. Soc. 
London Ser.~A} {\bf 455} (1999), 219-225.

\bibitem{edm2} D. E. Edmunds and J. R\'akosn\'{\i}k, Density of smooth functions in 
$W^{k,p(x)}(\Omega)$, {\it Proc. Roy. Soc. London Ser.~A} {\bf 437} (1992), 229-236.

\bibitem{edm3} D. E. Edmunds and J. R\'akosn\'{\i}k, Sobolev embedding with variable exponent, 
{\it Studia Math.} {\bf 143} (2000), 267-293.

\bibitem{FH} X. L. Fan and X. Han, Existence and multiplicity of solutions for $p(x)$-Laplacian 
equations in $\RR^N$, {\it Nonlinear Anal.} {\bf 59}  (2004), 173-188.

\bibitem{FZZ} X. L. Fan, Q. H. Zhang, and D. Zhao, Eigenvalues of
$p(x)$-Laplacian Dirichlet problem, {\it J. Math. Anal. Appl.}
{\bf 302} (2005), 306-317.

\bibitem{FZ1} X. L. Fan and D. Zhao, On the spaces $L^{p(x)}(\Omega)$ 
and $W^{m,p(x)}(\Omega)$, {\it J. Math. Anal. Appl.}, {\bf 263}
(2001), 424-446.

\bibitem{garper} J. P. Garcia Azorero and I. Peral Alonso, Existence and nonuniqueness for the 
$p$--Laplacian nonlinear eigenvalues,
{\it Comm. Partial Differential Equations} {\bf 12} (1987), 1389-1403.

\bibitem{jabri} Y. Jabri, {\it The Mountain Pass Theorem. Variants, Generalizations and Some 
Applications}, Encyclopedia of Mathematics and its Applications, Vol. 95, Cambridge University 
Press, Cambridge, 2003.

\bibitem{KR} O. Kov\'a\v cik and J. R\'akosn\'{\i}k, On spaces $L^{p(x)}$ and
$W^{1,p(x)}$, {\it Czechoslovak Math. J.} {\bf 41} (1991), 592-618.

\bibitem{lew} M. Lewin, A mountain pass for reacting molecules, {\it Ann. Henri Poincar\'e} {\bf 
5} (2004), 477-521.

\bibitem{M} J. Musielak, {\it Orlicz Spaces and Modular Spaces},
Lecture Notes in Mathematics, Vol. 1034, Springer, Berlin, 1983.

\bibitem{nak} H. Nakano, {\it Modulared Semi-ordered Linear Spaces}, Maruzen Co., Ltd., Tokyo, 
1950.

\bibitem{orl} W. Orlicz,  \"Uber konjugierte Exponentenfolgen, {\it Studia Math.} {\bf 3} (1931), 
200-212.

\bibitem{palsma} R. S. Palais and S. Smale, A generalized Morse theory, {\it Bull. Amer. Math. 
Soc.} {\bf 70} (1964), 165-171.

\bibitem{pre} R. Precup, An inequality which arises in the absence of the mountain pass geometry, 
{\it J.~Inequal. Pure Appl. Math. (JIPAM)} {\bf 3} (2002), no. 3, Article 32, 10 pp. (electronic).

\bibitem{rab} P. H. Rabinowitz, {\it Minimax Methods in Critical Point Theory with Applications to 
Differential Equations}, CBMS Regional Conference Series in Mathematics, Vol. 65, American 
Mathematical Society, Providence, RI, 1986. 

\bibitem{ruz} M. Ruzicka, {\it Electrorheological Fluids Modeling and Mathematical Theory}, 
Springer-Verlag, Berlin, 2000.

\bibitem{sha} I. Sharapudinov, On the topology of the space $L^{p(t)}([0;1])$,
{\it Matem. Zametki} {\bf 26} (1978),  613-632.

\bibitem{tse} I. Tsenov, Generalization of the problem of best approximation 
of a function in the space $L^s$, {\it Uch. Zap. Dagestan Gos. Univ.} {\bf 7} (1961), 25-37.

}
\end{thebibliography}
\end{document}